\documentclass[11pt]{article}

\usepackage[letterpaper,width=135mm,height=203mm]{geometry}

\usepackage{graphicx}

\newtheorem{theorem}{Theorem}
\newtheorem{lemma}[theorem]{Lemma}
\newtheorem{claim}[theorem]{Claim}

\newenvironment{Proof}{\proofing}{\QED}
\newcommand{\QED}{\hspace{8mm}\mbox{\textsc{qed}}\smallskip}
\newcommand{\proofing}{\textsc{Proof.}}

\newcommand{\Tcite}[1]{\mbox{\normalfont\textrm{\cite{#1}}}}

\newcommand{\curlyG}{\mathcal{G}}
\newcommand{\curlyE}{\mathcal{E}}
\newcommand{\iso}{\mathop{\mathit{\iota}}}
\newcommand{\Phat}{\widehat{P}}
\newcommand{\curlyX}{\mathcal{X}}

\makeatletter
\long\def\@caption#1[#2]#3{\begingroup \@parboxrestore 
\if@minipage \@setminipage \fi \normalsize \sffamily \@makecaption {\csname fnum@#1\endcsname }{\ignorespaces #3}\par \endgroup}
\makeatother

\begin{document}

\begin{center}
\textbf{\large Disjoint Isolating Sets and Graphs with Maximum Isolation Number} \\[3mm]
Geoffrey Boyer and Wayne Goddard \\
School of Mathematical and Statistical Sciences, Clemson University
\end{center}

\begin{abstract}
An isolating set in a graph is a set $X$ of vertices such that every edge of the graph is incident with a vertex of $X$ or its neighborhood.
The isolation number of a graph, or equivalently the vertex-edge domination number, is the minimum number of 
 vertices in an isolating set. Caro and Hansberg, and independently \.{Z}yli\'{n}ski, showed that 
 the isolation number is at most one-third the order for every connected graph of order at least $6$. 
 We show that in fact all such graphs have three disjoint isolating sets. Further, using a family introduced by 
Lema\'{n}ska, Mora, and Souto-Salorio,
we determine all graphs with equality in the original bound.
\end{abstract}

\section{Introduction}

Recently, 
Caro and Hansberg~\cite{CH} introduced the concept of isolation in a graph~$G$. 
Specifically, they considered a set $X$ of vertices such that the graph $G-N[X]$ 
has some sparseness property, where $N[X]$ denotes the closed neighborhood of $X$ (it together with its neighbors).
For example, if $G-N[X]$ has no vertex, then $X$ is a dominating set. The idea of
``paying'' for a vertex set but removing both it and its neighbors also appears in the vulnerability
literature; see for example~\cite{CW}.

As a special case of the general concept, Caro and Hansberg~\cite{CH} defined
an \textit{isolating set} of $G$ as a set $X$ of vertices such that $G-N[X]$ has no edge, and the
\textit{isolation number} of $G$, denoted $\iso(G)$, as the minimum cardinality of an 
isolating set. It should be noted that this parameter is equivalent to the \textit{vertex-edge domination number}
of $G$, as introduced earlier by Lewis et al.~\cite{LHHF}. They defined
a vertex~$u$ to \textit{$ve$-dominate} an edge~$f$ if either $f$ is incident to $u$,
or at least one of the ends of~$f$ is adjacent to $u$. Further, a set $S$ is
\textit{vertex-edge dominating} if every edge is $ve$-dominated by at least one vertex in $S$;
in other words, all edges of $G$ are incident to a vertex in~$N[S]$, or equivalently, that $S$ is an isolating set.
The isolation number is also the same as the $2$-distance vertex covering number, as defined by Canales et al.~\cite{CHMM}.

Caro and Hansbeg proved a bound on the isolation 
number of a graph, later also proved by \.{Z}yli\'{n}ski for vertex-edge domination number:

\begin{theorem} \Tcite{CH, Zylinski} \label{t:bound}
If graph $G$ is connected and not $K_2$ or $C_5$, then its isolation number is
at most one-third its order.
\end{theorem}

They and others posed the question of characterizing the graphs with isolation or vertex-edge domination number 
equal to one-third their order.
Earlier Krishnakumari et al.~\cite{KVK} had characterized trees of order $n$ with vertex-edge domination number equal to $n/3$, 
later also proved by Dapena et al.~\cite{DLSV} for isolation number.
Recently, Lema\'{n}ska et al.~\cite{LMS} introduced an infinite family of graphs 
$\curlyG$ such that every graph in the family has 
isolation number one-third its order. 
They showed that every unicyclic or block graph with this property
is in this family, apart from two exceptions. They then 
suggested that apart from this family there is only a finite number of connected graphs in general  with the property.
This we prove in Section~\ref{s:extremal}.

In another direction, it is well-known that every connected nontrivial graph has two disjoint dominating sets
but not necessarily three. We show in Section~\ref{s:disjoint} that every connected graph, except $K_2$ and $C_5$,
has three disjoint sets such that the vertices not dominated by all three sets form 
an independent set. 
As a corollary this means that 
for a connected graph, except for $K_2$ and~$C_5$, 
there are three disjoint (but not necessarily minimum)  isolating sets. 
This strengthens a result of Kiser~\cite{Kiser} (Theorem~4.8) who
showed that every connected graph of order at least~$3$
has three disjoint distance-$2$ dominating sets. A \textit{distance-$2$ dominating set}
is a set~$X$ such that every vertex of the graph is either in $N[X]$ or has a neighbor in $N[X]$; this is a weaker
property than isolating set.  

\newpage

\section{Three Disjoint Sets that Nearly Dominate} \label{s:disjoint}

We prove:

\begin{theorem} \label{t:disjoint}
If $G$ is a connected graph of order at least $3$ other than $C_5$, then 
there exists a partition $(A_1,A_2,A_3)$ of $V(G)$ such that 
$\curlyX = \bigcup_{i=1}^3 \left( V(G) - N[A_i] \right)$ is an independent set.
\end{theorem}

Note that the conclusion is slightly stronger than just $V(G) - N[A_i]$ being an independent set for each $i$.

We refer to the partition as a \textit{coloring}.
The proof of the theorem is by induction on the order and is divided into a sequence of claims.
Furthermore, for the induction, note that the conclusion is also true for $G=K_1$ if one relaxes
partition to weak partition.
The base case is included in the following claim.

\begin{claim}
If $G$ is a star with at least three vertices, or 
a cycle $C_n$ with $n\neq 5$, then there exists such a coloring.
\end{claim}
\begin{Proof}
For the star, any coloring using all three colors is valid. So assume $G$ is $C_n$ and 
let $r$ be the remainder when $n$ is divided by $3$.
If $r=0$ then use the coloring that repeats $1,2,3$;
if $r=1$ then use the coloring that is $1,2,1,3$ followed by repeating $1,2,3$; and
if $r=2$ then use the coloring that is $1,3,2,1,3$ followed by repeating $1,2,3$.
If $r=0$ then $\curlyX$ is empty; otherwise $\curlyX$ contains two vertices and these are at distance at least two on the cycle.
\end{Proof}

So assume $G$ is neither a star nor a cycle.
 
\begin{claim} \label{c:cycleInduct}
If $G$ has a cycle $C$ of length a multiple of $3$, then one can induct.
\end{claim}
\begin{Proof}
Color the cycle $C$ with the $1,2,3$ coloring as above. It follows that no vertex of $C$ is in $\curlyX$.
Let $F$ be a component of $G-C$. If $F$ has a valid coloring then use it.
If $F$ is $K_2$, say $x_1x_2$ with $x_1$ having a neighbor $v$ on $C$,
color $x_1$ and $x_2$ so that $x_1$, $x_2$, and $v$ have different colors. Only $x_2$ can be in~$\curlyX$.
If $F$ is $C_5$, say $x_1\ldots x_5$ with $x_1$ having a neighbor $v$ on $C$, 
then give $x_3$ the same color as $v$, give $x_2$ and $x_5$ the second color, and give $x_1$ and $x_4$ the third color.
Only vertex $x_5$ can be in $\curlyX$.
\end{Proof}

This idea can be re-used as follows:

\begin{claim} \label{c:path}
If $G$ has a nontrivial path $P = v_1 \ldots v_k$ such that $G-P$ is disconnected,
and there exists a neighbor $x$ of $v_1$ and a neighbor $y$ of $v_k$ in different components of $G-P$,
then one can induct.
\end{claim}
\begin{Proof}
Color the path with the $1,2,3$ coloring. All vertices of $P$ except $v_1$ and $v_k$ are dominated
by all three colors.
Then color the remainder of the graph, as in the previous claim, choosing the color of 
$x$ to be $3$ and the color of $y$ to be the color missing at $v_k$; hence no vertex of $P$ is in $\curlyX$.
\end{Proof}

\begin{claim}
We may assume $G$ is $2$-connected and has minimum degree at least $3$.
\end{claim}
\begin{Proof}
Suppose $G$ has a cut-vertex $v$. Then since $G$ is not a star, 
vertex $v$ has a neighbor $z$ in some nontrivial component of $G-v$.
Then $vz$ is a separating path, and one can apply the above claim, to induct.
So we may assume $G$ is $2$-connected.

Suppose $G$ has a vertex of degree $2$. Consider the shortest cycle $C$
containing such a vertex; necessarily induced. Since $G$ is not a cycle, the cycle contains a vertex
with a neighbor outside the cycle. Thus there must exist consecutive vertices
$v$ and~$w$ on $C$ such that $v$ has degree $2$ in $G$ and $w$ has a neighbor outside the cycle.
Let $P$ be the path $C-v$; this satisfies the hypothesis of the above claim, and so we can induct.
\end{Proof}

\begin{claim}
If $G$ has a cycle $C$ such that $G-C$ is disconnected, then one can induct. 
\end{claim}
\begin{Proof}
If every vertex of $C$ has a neighbor outside $C$, then there must be two consecutive vertices on $C$
with neighbors in different components of $G-C$. Thus we can apply Claim~\ref{c:path}.
Otherwise, assume there is a vertex on $C$ that has no neighbor outside $C$. Then choose such a vertex $v$
such that one of its neighbors on $C$ does have a neighbor outside $C$. The path $C-v$ satisfies
the conditions of Claim~\ref{c:path} since $v$ is isolated by the path's removal.
\end{Proof}

Thomassen determined the graphs that have no separating cycle. In particular, Theorem 3.1a of \cite{Thomassen}
can be summarized as:

\begin{theorem} \Tcite{Thomassen} \label{t:noSeparating}
If $G$ is a $2$-connected graph with minimum degree at least $3$ and with no cycle $C$ such that $G-C$ is disconnected, then 
$G$ is one of: the complete graph $K_n$, the wheel~$W_n$, the balanced or almost balanced complete bipartite graph 
possibly with some specific edges added, the prism on six vertices, or one of the $23$ graphs given in Figure~1 of that paper.
\end{theorem}

It can then be checked that all the graphs given in Theorem~\ref{t:noSeparating} have either a $3$-cycle or a $6$-cycle, and so allow induction
by Claim~\ref{c:cycleInduct}. 
This completes the proof of Theorem~\ref{t:disjoint}.

\section{Graphs with Maximum Isolation Number} \label{s:extremal}

\begin{figure}[h]
\begin{center}
\includegraphics{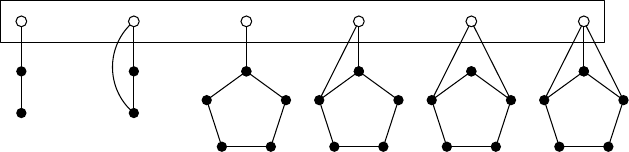}
\end{center}
\caption{A graph in $\curlyG$} 
\label{f:curlyG}
\end{figure}

A graph in the family $\curlyG$ as given by Lema\'{n}ska et al.~\cite{LMS} is defined as follows. There is a connected \textit{base} graph. Associated
with each vertex of the base graph, which we call a \textit{hook}, one introduces one pendant, which
is either $K_2$ or $C_5$.  For a $K_2$-pendant, the hook is joined to one or two vertices.
For a $C_5$-pendant, the hook is joined to one, or any two, or three consecutive vertices of the $5$-cycle
(that is, any nonempty set that is not a vertex cover). See Figure~\ref{f:curlyG}.

Lema\'{n}ska et al.~\cite{LMS} observed that these graphs have isolation number one-third their order. 
It is also immediate that something slightly stronger is true:

\begin{lemma}
For a graph of $\curlyG$ of order $n$, both its
domination number and isolation number are $n/3$.
\end{lemma}

We prove:

\begin{theorem} \label{t:mainSummary}
For $n\ge 15$, all connected graphs of order $n$ with isolation number~$n/3$ are in $\curlyG$.
\end{theorem}

\begin{figure}[h]
\begin{center}
\includegraphics{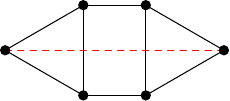}
\quad
\includegraphics{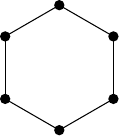}
\bigskip

\includegraphics{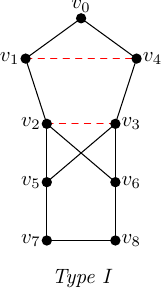}
\quad
\includegraphics{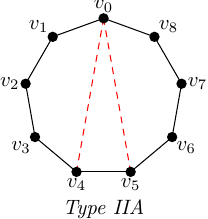}
\quad
\includegraphics{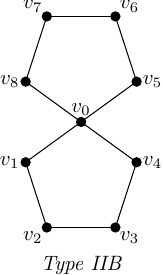}
\bigskip

\includegraphics{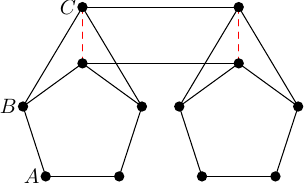}
\end{center}
\caption{The graphs in $\curlyE$} 
\label{f:curlyE}
\end{figure}

For the proof it is helpful to determine all examples of equality.
We define a family $\curlyE$ of $14$ exceptional graphs, shown in Figure~\ref{f:curlyE}.
Note that the dashes represent optional edges, and thus some depictions represent multiple graphs.
Up to isomorphism there are three graphs of order $6$,
eight graphs of order $9$, and three graphs of order $12$. (It is perhaps interesting to note that these
graphs also have domination number equal to isolation number.) 

We prove the full characterization:

\begin{theorem} \label{t:main}
A connected graph has isolation number one-third its order if and only if it is
in $\curlyG \cup \curlyE$. 
\end{theorem}

The proof is by induction on the order. 

\subsection{Preparing for Induction}

The base case is orders $3$, $6$, and $9$.
The order-$3$ case is trivial and the order-$6$ case was noted in~\cite{LMS}.

\begin{lemma}
The connected graphs of order $9$ with isolation number $3$ are the $18$ graphs in $\curlyG$ and the $8$ graphs in
$\curlyE$.
\end{lemma}
\begin{Proof}
By computer search. We enumerate all connected graphs of order $9$ and determine which do not have an isolating set of cardinality $2$.
Java code is provided on \textsf{https://people.computing.clemson.edu/$\sim$goddard/papers/isolation/}.
\end{Proof}

For the inductive step, it is helpful to have a structural result.
(This is somewhat similar to the idea used in the induction proof in~\cite{CH}.) 

\begin{lemma} \label{l:structure}
If $G$ is a connected graph of order at least $3$, then there exists a not necessarily induced star $S$ 
on at least $3$ vertices such that $G-S$ has at most one nontrivial component.
\end{lemma}
\begin{Proof}
Consider a spanning tree $T$ of $G$. If $T$ is a star, we are done; so assume $T$ has diameter at least $3$.
Root $T$ at the beginning of a longest path; 
say it starts at vertex $a$ and ends with $xyz$.

If $y$ has two or more children, then let $S_1$ be the star of $y$ 
and all its children; since $T-S_1$ is connected it follows that $G-S_1$ is connected.
So assume $y$ has only one child. 
Note that we may similarly assume all children of $x$ have degree either one or two in $T$.

If $x$ has degree $2$ in $T$, then consider the star $S_2$ of $x,y,z$;
again $G-S_2$ is connected.
So assume $x$ has a second child. If there is an edge in $G$ between
grandchildren of $x$, say $zz'$, then $G-S_3$ is connected where $S_3$ is the star of $z$, $z'$, and $y$.

So assume there is no edge between grandchildren of $x$. Consider the star $S_4$ of $x$ and all its children. Its removal from $G$ leaves 
a component $C$ containing vertex $a$, and potentially components containing
grandchildren of $x$. But since there is no edge between grandchildren, each grandchild is either in a 
component by itself, or is joined by an edge to~$C$. So $S_4$ is the desired star.
\end{Proof}

So assume 
$G$ is a graph of order $n\ge 12$ with $\iso(G)=n/3$.
By Lemma~\ref{l:structure}, there is a star $S$ with at least three vertices such that $H=G-S$ 
has one nontrivial component and possibly some isolates.
Say the star $S$ is centered at $x$ with neighbors $y_1, \ldots, y_d$. 
Then one can add $x$ to the isolating set of $H$ and so $\iso(G) \le \iso(H)+1$. 
Since $H$ has $n-d-1$ vertices,  
by Theorem~\ref{t:bound} applied to the nontrivial component of $H$, 
it follows that $d=2$, $H$ has no isolate, and $\iso(H) = (n-3)/3$.
By the induction hypothesis, the graph $H$ is in $\curlyE \cup \curlyG$.

\subsection{When $H$ is in $\curlyE$}

Assume $H$ is one of the exceptional graphs. We will need the following two ``extendability'' lemmas:

\begin{claim} \label{c:12extend}
If $H$ is one of the exceptional graphs of order $12$, then any two vertices of $H$ can be extended to 
an isolating set of $H$ of cardinality $4$.
\end{claim}
\begin{Proof}
Say vertices $z_1$ and $z_2$ are given.
Call the degree-$2$ vertices that have a degree-$2$ neighbor \textit{$A$-vertices},
the degree-$3$ vertices with an $A$-neighbor the \textit{$B$-vertices}, and the remaining four vertices the \textit{$C$-vertices}. 
We use ``half'' to mean one of the two components that would result if the edges joining $C$-vertices were to be removed.
If $z_1$ and $z_2$ are from different halves,
then let $b_1$ be a $B$-vertex distinct from $z_1$ and $b_2$ a $B$-vertex  distinct from $z_2$;
it can be checked that $\{z_1,z_2,b_1,b_2\}$ is an isolating set of $H$.
If $z_1$ and $z_2$ are in the same half, there are two cases. If they are the two $C$-vertices,
then add an $A$-vertex from each half. Otherwise, add a $B$-vertex and $C$-vertex from the other half,
choosing a $C$-vertex whose neighbor in the first half is not already dominated if there is such a one. 
In both cases the quartet is an isolating set of $H$.
\end{Proof}

In the following, we use the labels and types as shown in Figure~\ref{f:curlyE}.

\begin{claim} \label{c:9extend}
If $H$ is one of the exceptional graphs of order $9$, then any two vertices of $H$ can be extended to 
an isolating set of $H$ of cardinality $3$ except possibly pairs $\{v_0,v_1\}$, $\{v_0,v_4\}$, $\{v_7,v_8\}$ for type~I,
and  pairs $\{v_2,v_3\}$, $\{v_6,v_7\}$ for type II.
\end{claim}
\begin{Proof}
The following tables provide one possible third vertex for each pair, where in the second table, the notation
$a/b$ means vertex $v_a$ works for type IIA and vertex $v_b$ for type~IIB.

{\small \begin{tabular}{c|cccccccc}
\multicolumn{9}{c}{Type I} \\
 &  1 & 2 & 3 & 4 & 5 & 6 & 7 & 8\\ \hline
0 &   --  & 5 & 5 & --  & 2 & 2 & 2 & 2 \\
1 &    & 5 & 5 & 7 & 2 & 2 & 3 & 3 \\
2 &    &  & 5 & 5 & 0 & 0 & 0 & 0 \\
3 &    &  &  & 5 & 0 & 0 & 0 & 0 \\
4 &    &  &  &  & 2 & 2 & 1 & 1 \\
5 &    &  &  &  &  & 0 & 0 & 0 \\
6 &    &  &  &  &  &  & 0 & 0 \\
7 &    &  &  &  &  &  &  & --  \\
\end{tabular}
\qquad
\begin{tabular}{c|cccccccc}
\multicolumn{9}{c}{Type II} \\
 &  1 & 2 & 3 & 4 & 5 & 6 & 7 & 8\\ \hline
0   & 5 & 5 & 5 & 5 & 1 & 2 & 3 & 4 \\
1   &  & 6 & 6 & 6 & 0 & 2 & 3 & 4/0  \\
2   &  &  & --  & 7 & 0 & 0 & 4 & 5 \\
3   &  &  &  & 8/7  & 0 & 0 & 0 & 5 \\
4   &  &  &  &  & 0 & 0 & 0 & 0 \\
5   &  &  &  &  &  & 1/2  & 2 & 2 \\
6   &  &  &  &  &  &  & --  & 2 \\
7   &  &  &  &  &  &  &  & 3 \\
\end{tabular}} \\
\end{Proof}

We continue with the overall proof. Recall that $S$ has vertex set $\{x,y_1,y_2\}$ with $x$ adjacent to $y_1$ and $y_2$.

\begin{claim}
Vertex $x$ has no neighbor in $H$.
\end{claim}
\begin{Proof}
Suppose $x$ has neighbor $z$ in $H$.
The vertex $x$ can be added to any isolating set of $H-z$ to form an isolating set of $G$. 
Note that for all but one graph, an exceptional $H$ is $2$-connected.
In that case, $H-z$ is connected,
so its isolation number is at most $(n-4)/3$, and thus $G$ is not extremal.
A similar claim holds for the remaining possibility except when $z$ is the cut-vertex thereof.
But in that case $H-z$ has no component of order $2$ or $5$, and so again its isolation number is at most $(n-4)/3$.
\end{Proof}

\begin{claim}
Both $y_1$ and $y_2$ have neighbors in $H$.
\end{claim}
\begin{Proof}
Suppose that $y_1$ say has no neighbor in $H$. Then $y_2$ must have a neighbor in $H$, say $z$.
Again $H-z$ has isolation number at most $(n-4)/3$, and an isolating set of $H-z$ can be extended 
to one of $G$ by adding $y_2$, so that $G$ is not extremal.
\end{Proof}

So assume $y_1$ has neighbor $z_1$ and $y_2$ has neighbor $z_2$ in $H$ (possibly $z_1=z_2$).
Consider first the case that $H$ has order $12$.
Then by Claim~\ref{c:12extend}, one can extend $\{z_1, z_2\}$ to an isolating set $J$ of $H$ 
of size $4$. The resultant $J$ is also an isolating set of $G$, since only $x$ remains of the star $S$; 
and thus $G$ is not extremal.

Consider the case that $H$ has order $9$.
Again, if $\{z_1,z_2\}$ can
be extended to an isolating set $J$ of $H$ of cardinality $3$, then the set $J$ is an isolating
set of $G$ and hence $G$ is not extremal. 

Consider $H$ of type II. Since by Claim~\ref{c:9extend} the two possible ``bad'' pairs are disjoint, the only possibility is that both $y_1$ and $y_2$ have 
unique neighbors in $H$. But if $y_1$ has neighbor $v_2$ and $y_2$ has neighbor $v_3$, then $\{v_1,v_6,y_2\}$ is an isolating set of $G$;
and if $y_1$ has neighbor~$v_6$ and $y_2$ has neighbor $v_7$, then $\{v_1,v_5,y_2\}$ is an isolating set of $G$. 
That is, the resultant $G$ is not extremal.
Consider $H$ of type~I.  If both $y_1$ and $y_2$ have unique neighbors in $H$,
then again it can be checked that the resultant $G$ is not extremal.
But if $y_1$ has neighbor~$v_0$ and $y_2$ has neighbors $v_1$ and $v_4$, then for each $H$
one obtains an extremal graph. However, note that the $G$ with edge $v_1v_4$ present
but not $v_2v_3$ is isomorphic to the $G$ with edge $v_2v_3$ present but not $v_1v_4$. So one obtains 
that $G$ is one of
the three graphs of order~$12$ in $\curlyE$.


\subsection{When $H$ is in $\curlyG$}

Assume $H$ is in $\curlyG$.
Let $B$ be the vertices of the base graph of $H$, and let $C$ be a minimum isolating set of $H$ that contains $B$. (It can easily be checked that one
can choose such a minimum isolating set). For a pendant $P$, we use the notation $\Phat$ for the set of its vertices and its hook.

\begin{claim} \label{c:xNoPendant}
Vertex $x$ has no neighbor in the pendants of $H$.
\end{claim}
\begin{Proof}
Suppose $x$ has neighbor $z$ in the pendants.
Then vertex $x$ can be added to any isolating set of $H-z$ to form one of $G$. 
Since $H-z$ has no component of order $2$ or $5$,
its isolation number is at most $(n-4)/3$, and thus $G$ is not extremal.
\end{Proof}

\begin{claim}
We may assume vertices $y_1$ and $y_2$ are not adjacent.
\end{claim}
\begin{Proof}
Suppose $y_1y_2$ is an edge. Then, since one can choose $x$ as any of the trio, it follows 
from Claim~\ref{c:xNoPendant}
that 
none of them has neighbors in the pendants. That is, their only neighbors outside $S$ are in $B$.
If two of them have neighbors in $B$, then $C$ is an isolating set of $G$, a contradiction.
Thus only one of the trio has a neighbor in $B$, and so $G$ is in $\curlyG$.
\end{Proof}

\begin{claim} \label{c:pendantNeigh}
We may assume both $y_1$ and $y_2$ have a neighbor in the pendants of $H$.
\end{claim}
\begin{Proof}
Suppose $y_2$ say has no neighbor in the pendants. There are two cases:

\textit{Case 1:} Assume $y_2$ has degree $1$. \\
If all of $y_1$'s neighbors are in $B$, then since $C$ is not an isolating set of $G$, it follows that $x$ has no neighbor in $B$; thus
$G$ is in $\curlyG$. So assume
$y_1$ has a neighbor in a pendant $P$ with hook vertex $v$. 
Let $C' = C - \Phat$ and consider deleting $N[C' \cup \{y_1\} ]$ from $G$.
Then only $y_2$ remains of $S$, and so all vertices outside $P$ that remain are guaranteed to be isolated. 
Further, vertex $v$ is adjacent to $y_1$ and so does not remain.
If $P$ was a $K_2$-pendant, then one isolated vertex of $P$ remains; if $P$ was a $C_5$-pendant, then what remains has a vertex $a$ of degree $2$ 
and removing $N[a]$ leaves no edges. 
That is, $G$ has an isolating set of size $n/3-1$ and so is not an extremal graph, a contradiction. 

\textit{Case 2:} Assume $y_2$ has a neighbor in $B$. \\
If $y_1$ has a neighbor in $B$, then the set $C$ would be an isolating set of~$G$, a contradiction.
So we may assume every neighbor of $y_1$ except $x$ is in the pendants of $H$. 

Case 2a) Assume $y_1$ has neighbors $z_1$ and $z_2$ in different pendants. Then deleting $\{ y_1, z_1, z_2, x \}$ from $G$
yields a graph $G'$ with no $2$- or $5$-component except possibly the component
containing $B$ and $y_2$. But it is easily checked that if that component has order $5$ 
then $|B|=2$ and both pendants are $K_2$-pendants, which is impossible since $|H|\ge 9$.
That is, $\iso(G) \le (n-4)/3+1$, since $y_1$ can be added to any isolating set of $G'$,
and so $G$ is not extremal.

Case 2b) Assume $y_1$ has neighbor(s) only in one pendant, say $P$ with hook~$v$. Let $C' = C - \Phat$. 
Then removing $N[C'] - S$ from $G$ yields isolates except for the component $F$ induced by
$V(P) \cup V(S)$.  If $P$ is a $C_5$-pendant, then $F$ has order~$8$ and hence has an 
isolating set of two vertices. If $P$ is a $K_2$-pendant, then $F$ has order $5$; but $y_2$ has degree $1$ in $F$ and so $F$ is
not a $5$-cycle and hence has a $1$-element isolating set. It follows that $G$ is not extremal, a contradiction.
\end{Proof}

\begin{claim} \label{c:notJustOne}
If there exists one pendant $P$ containing all neighbors of $y_1$ and $y_2$ in the pendants of $H$, then we are done.
\end{claim}
\begin{Proof}
Let $C' = C - \Phat $.
Assume first that $y_1$ and $y_2$ have a common neighbor say $z$. Then $C'$ can be extended to 
an isolating set $C''$ of $H$ and hence $G$ by adding~$z$ if 
$P$ is $K_2$, or by adding $z$ and a vertex of $P$ otherwise; in either case $|C''|=|C|=(n-2)/3$, a contradiction.

Assume second that $y_1$ and $y_2$ have no common neighbor in $P$.
Say $z_1$ is a neighbor of $y_1$ and $z_2$ is a neighbor of $y_2$.
Then adding $z_1$ and $z_2$ to $C'$ yields an isolating set of $H$ and hence of $G$; so we are done unless
$P$ was a $K_2$. But in that case, one can easily check that since $C$ is not an isolating set of $G$, adding the vertices of $S$ to $P$ turns it into a $C_5$-pendant.
Thus $G$ is in $\curlyG$.
\end{Proof}

By Claims~\ref{c:pendantNeigh} and~\ref{c:notJustOne}, we may assume that
$y_1$ has a neighbor in some pendant~$P_1$ with hook~$v_1$, and 
$y_2$ has a neighbor in some pendant $P_2$ with hook~$v_2$.

\begin{claim}
If $|B|>2$, then we are done.
\end{claim} 
\begin{Proof}
Assume first that there is another pendant of $H$ containing a neighbor of $y_1$ or~$y_2$.
Say $z_3$  is a neighbor of $y_1$ in pendant $P_3$.
Then consider the graph $G' = G - \{ y_1, z_1, z_3, x \}$.
We claim that $G'$ has no component of order $2$ or $5$. To see this, observe that no
vertex is removed from the base, and $y_2$ is in the component of $G'$ containing the base.
Removing $z_1$ or $z_3$ can separate a component in $P_1$ or $P_3$, but such a component has order
$1$ or $4$. It follows by Theorem~\ref{t:bound} that the isolation
number of $G'$ is at most $(n-4)/3$, and so $G$ is not extremal.

So assume neither $y_1$ nor $y_2$ has a neighbor in any other pendant.
Let $C' = C  - ( \Phat_1 \cup \Phat_2)$.
Then, since $|B|>2$, the set $C'$ is not empty.
If one removes from $H-( \Phat_1 \cup \Phat_2)$ the vertices of $N[C']$ therein, then all remaining vertices are from pendants, and they are 
now isolated. Further, none of the remaining vertices has a neighbor in $\Phat_1$ or~$\Phat_2$, by the construction of graphs in~$\curlyG$;
and none has a neighbor in $S$ by the above discussion. 
That is, if one removes $N[C']$ from $G$, all remaining vertices outside $\Phat_1 \cup \Phat_2 \cup S $ are isolated.
Further, since the base graph is connected, at least one of $v_1$ and $v_2$ is removed. 
Because $S$ joins the two pendants, what remains of $\Phat_1 \cup \Phat_2 \cup S$ is connected, and its order is not a multiple of $3$. 
It follows that $G-N[C']$ has
isolation number less than one-third its order, which implies that $G$ is not extremal.
\end{Proof}

\begin{claim}
If $|B|=2$, then we are done.
\end{claim} 
\begin{Proof}
Since $H$ has order at least $9$, there are two cases.

\textit{Case 1:} Assume both $P_1$ and $P_2$ are $C_5$-pendants.\\
Then there is a vertex $a_1$ in $P_1$ such that
$\{z_1,a_1\}$ dominates all of $\Phat_1$ except possibly some vertex $c_1$ that has no neighbor in $P_2$.
Similarly there is a vertex $a_2$ in $P_2$ such that $\{z_2,a_2\}$ dominates 
all of $\Phat_2$ except possibly one vertex $c_2$ on the cycle. But 
the vertices $c_1,c_2,x$ are independent. (Recall that by 
Claim~\ref{c:xNoPendant} vertex $x$ has no neighbor in the pendants.)
Thus, the quartet $\{z_1, z_2, a_1, a_2 \}$ is an
isolating set of $G$ of cardinality $4$, a contradiction.


\textit{Case 2:}  Assume $P_1$ is a $C_5$-pendant and $P_2$ a $K_2$-pendant.\\
Then by the same reasoning as before, there is a vertex $a_1$ in $P_1 $ so that
$\{z_1,a_1\}$ dominates all of $\Phat_1$ except possibly some vertex $c_1$, and $c_1$ has no neighbor in~$P_2$.
If one removes $N[ \{ z_1, a_1, z_2 \}]$ from $G$, the only vertices that can remain are $c_1$, $x$, and $v_2$. 
Since the isolation number of $G$ is $4$, that removal 
must leave an edge. The only possibility for such an edge is $xv_2$; it follows thence that $v_2$ is not adjacent to $z_2$.
See Figure~\ref{f:portion}.

\begin{figure}[h]
\begin{center}
\includegraphics{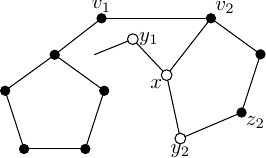}
\end{center}
\caption{A portion of an extremal graph} 
\label{f:portion}
\end{figure}

Consider now the star $S' = P_2 \cup y_2$. Since $y_1$ has an edge to $P_1$, the graph
$H'=G-S'$ is bridgeless and connected. So if $H'$ is an extremal graph, it is in~$\curlyE$,
and that case was already handled.
\end{Proof}

This concludes the proof of Theorem~\ref{t:main}.

\section{Further Thoughts}

We note that all graphs in $\curlyG$ with nontrivial base graph have cut-vertices. So this raises the question of what is the asymptotic upper bound on the isolation number 
in $2$-connected graphs. Another direction is to obtain better bounds for specific families of graphs, such as regular graphs (see~\cite{ZZ}).
Further, the relationship between domination and isolation deserves further investigation. For example, the corona of the complete
graph has maximum domination number but minimum isolation number.


\end{document}